

\input amstex
\documentstyle{amsppt}
\magnification\magstep1
\NoRunningHeads
\hsize6.9truein
\pageheight{21 truecm}
\baselineskip=15pt

\catcode`@=11
\def\logo@{\relax}
\catcode`@=\active

\def\du{\coprod}

\def\ol#1{\overline{#1}}
\def\Spec{\operatorname{Spec}}
\def\ring#1.{\Cal O_{#1}}
\def\Hilb{\operatorname{Hilb}}
\def\Hom{\operatorname{Hom}}
\def\xdim{\operatorname{dim}}
\NoBlackBoxes
\def\subscheme{\subset}
\def\xk{\bold k}

\def\gc{GC }
\topmatter
\title Quotients by Groupoids \endtitle
\author Se\'an Keel and Shigefumi Mori \endauthor
\date Aug. 25, 1995 \enddate
\endtopmatter
\heading \S 1 Introduction and Statement of Results
\endheading

The main result of this note is the following:
(all relevant terms will be defined shortly.
The spaces are algebraic spaces of finite type over
a locally Noetherian base, unless otherwise mentioned.)

\proclaim{1.1 Theorem} Let $R \rightrightarrows X$
(i.e. $j : R \to X \times X$) be a flat groupoid
such that its stabilizer
$j^{-1}(\Delta_X) \rightarrow X$  is finite. Then there is
an algebraic space which is a uniform geometric quotient,
and a uniform categorical quotient. If $j$ is finite,
this quotient space is separated.
\endproclaim

Our principal interest in (1.1) is:
\proclaim{1.2 Corollary} Let $G$ be a flat group scheme acting
properly on $X$ with finite stabilizer, then a uniform geometric
and uniform categorical quotient for $X/G$
exists as a separated algebraic space.
\endproclaim

A version of (1.2) with stronger assumptions is the
main result of \cite{Koll\'ar95}, which we refer to for
applications, including  the existence of many moduli spaces
whose existence was previously known only in charateristic zero.
Our proof is based on entirely different ideas. We make use
of the greater flexibility
of groupoids vs. group actions to give an elementary and
straightforward argument.

Note that an immediate corollary of (1.1) is:

\proclaim{1.3 Corollary} (1) A separated algebraic stack
\cite{FaltingsChai80, 4.9} has a coarse
moduli space which is a separated algebraic space.

(2) An algebraic stack in Artin's sense \cite{Artin74, 5.1, 6.1}
has a GC quotient as an algebraic space if its
stabilizer is finite (2.7). If the stack is separated, so
is the GC quotient.
\endproclaim

(1.3.1) has a sort of folk status. It appears for example
in \cite{FaltingsChai80, 4.10}, but without reference or proof.

We note in \S 9 that for GIT quotients, the assumptions of
(1.1) are satisfied, and the quotients obtained are the same.

Several special cases of (1.1) appear in \cite{SGA3,V.4}, and
our (5.1) can be derived from \cite{SGA3, V.4.1}.

\definition{1.4 Definitions and Notations}
When we say a scheme in this paper,
it is not necessarily separated.
By a sheaf, we mean a sheaf in the qff (quasi-finite
flat) topology of schemes.

We fix a base scheme $L$ which is
locally Noetherian, and work only on $L$-schemes.
Thus the product $X \times Y$ means $X \times_L Y$, and if
we talk about the properties of $X$ (e.g. separated, of
finite type etc.), it is about those of $X \to L$ unless
otherwise mentioned.

By a geometric point, we mean
the spectrum of an algebraically closed field which is an
$L$-scheme. (This is not of finite type.)
\enddefinition

\definition{1.5 Definition} By a {\bf relation} we mean any map
$j: R \rightarrow X \times X$. We say $j$ is a {\bf pre-equivalence
relation} if the image of
$j(T) : R(T) \rightarrow X(T) \times X(T)$
is an equivalence relation (of sets) for all schemes $T$. If
in addition $j(T)$ is always a monomorphism, we call $j$ an
{\bf equivalence relation}. For a pre-equivalence
relation we write $X/R$ for the quotient sheaf.
\enddefinition

Throughout this paper $j$ will indicate a relation,
with projections $p_i: R \rightarrow X$.

\definition{1.6 Definition} A sheaf $Q$ is said to be
an {\bf algebraic space} over $L$ if $Q = U/V$ for some
schemes $U, V$ over $L$ and an equivalence relation
$j : V \to U \times U$ of finite type such that
each of the projections $p_i: V \rightarrow U$ is \'etale.
We say that $Q$ is {\bf of finite type} (resp. {\bf separated})
if furthermore $U, V, j$ are
chosen so that $U$ is of finite type
( resp. $j$ is a closed embedding).
\enddefinition

We note that fiber products always exist in the category of algebraic
spaces [Knudson71, II.1.5]. Therefore the notions of relation,
pre-equivalence relation and equivalence relation make sense when
$X,R$ are algebraic spaces.

\remark{1.7 Remark}
If $X/R \to Q$ is a map to an algebraic
space and if $Q' \to Q$ is a map of algebraic spaces, then
any relation $j : R \to X \times X$ can be pulled back to relation
$R'=R \times_Q Q'$ of $X'=X \times_Q Q'$
with $p_{i,R'} = p_{i,R} \times_Q Q' : R' \to X'$.
It is easy to check that if $R$ is a pre-equivalence relation
(or an equivalence relation) then so is $R'$. Another form of
pullback (which we call restriction) will be discussed in \S 2.
\endremark

\definition{1.8 Definition} Let $j:R \rightarrow X \times X$ be a
pre-equivalence relation. And let $q: X/R \rightarrow Q$ be a
map to an algebraic space $Q$. Consider the following properties
\roster
\item"(G)" $X(\xi)/R(\xi) \rightarrow Q(\xi)$ is a bijection for any
geometric point $\xi$.
\item"(C)" $q$ is universal for maps to (not necessarily separated)
algebraic spaces.
\item"(UC)" $(X \times_Q Q')/(R \times_Q Q') \rightarrow Q'$
satisfies (C) for any flat map $Q' \rightarrow Q$.
\item"(US)" $q$ is a universal submersion.
\item"(F)" The sequence of sheaves in the \'etale topology
$$0 \rightarrow \ring Q. \rightarrow q_*(\ring X.) \overset
{p_1^* - p_2^*} \to \rightarrow (q \circ p_i)_*(\ring R.)
$$
is exact (i.e. the regular functions of $Q$ are the
$R$-invariant functions of $X$).
\endroster

Note from the definition of quotient sheaf, that (UC) implies
(F).

If $q$ satisfies (C) it is called a {\bf categorical quotient}, and if
it satisfies (UC) it is called a {\bf uniform categorical quotient}.
If it satisfies (G) and (C) it is called a {\bf coarse moduli space}.
If it satisfies (G), (US) and (F) it is called a
{\bf geometric quotient}.
By a {\bf GC quotient} we will mean a quotient satisfying all
the above properties.
\enddefinition

\remark{1.9 Remark} Note (1.8.G) is by definition universal,
i.e. is preserved by any pullback $Q' \rightarrow Q$. Thus if
$U' \subset X' = X \times_Q Q'$ is an $R'$ invariant set
($R' = R \times_Q Q'$), then $U' = {q'}^{-1}(q'(U'))$ (as sets).
If the projections $p_i$ are universally open, and $q$ satisfies
(1.8.G) and (1.8.US), then $q$ is universally open.
\endremark

In the category of schemes a geometric quotient is
always categorical, and in particular unique
(see 0.1 of \cite{MumfordFogarty82}). This however fails for
algebraic spaces (see \cite{Koll\'ar95}).

Let us now give a rough idea of the proof, and an overview of
the layout of the paper:
The basic idea in our proof of 1.1 is to simplify the situation
by restriction. The idea comes from the sketch on pg. 218 of
\cite{MumfordFogarty82} of a proof of 1.1 in the analytic case, and
is roughly as follows:
To form the quotient we have
to identify elements in the same orbit. If $W$ is a general slice
through $x \in X$, then $W$ meets each orbit a finite number
of times, and we should be able to form the quotient from $W$ by
identifying each $w \in W$ with finitely many
``equivalent'' points (points in the same orbit in $X$). This
equivalence is no longer described by a group action, but
a good deal of structure is preserved. This leads to the
notion of a groupoid, defined in \S 2. The groupoid formalism
is very flexible, and allows us to simplify by various restrictions
(see \S 3). We can work \'etale locally around a point $x \in X$.
By taking the slice $W$ we reduce to the case
where the projections $p_1, p_2: R \rightarrow X$ are quasi-finite.
Then in \S 4 we reduce to the case where $R$ is {\it split}, a
disjoint union of a finite flat sub-groupoid $P$  and a piece
which does not have effects on $x$ (in particular $P$ contains
the stabilizer of $x$). In \S  5 we treat the case of
a finite flat groupoid, very similar to the case where a finite
group acts. In \S 7 we construct the quotient. We first mod out
by $P$, and afterwards we have a free \'etale action, which thus
defines an algebraic space.

Thanks: We would like to thank A.~Corti, J.~Koll\'ar, E.~Viehweg,
N.~Shepherd-Barron,
and P.~Deligne for helpful discussions. In addition the preprints
\cite{Viehweg95} and \cite{Koll\'ar95} helped clarify a number of
issues.
In particular, (3.1.2) comes from \cite{Viehweg95, 9.7}, and
definitions (2.2) and (2.4) come from \cite{Koll\'ar95, 2.11,2.14}.

\heading \S 2 Groupoids \endheading

We begin with some trivial categorical remarks, which in the end
is all there is to a groupoid.

First, for a small category $C$, write $R= Hom(C)$ and $X = Obj(C)$.
Then we have two natural maps $s,t: R \rightarrow X$ giving
for each morphism $f$ its source and target objects, and composition
gives a map $R \times_{(s, t)} R @>{c}>> R$, and
the identity map gives $e: X \rightarrow R$, a section for both
$s$ and $t$. There are obvious compatibilities between these
maps reflecting the axioms of a category, and a small
category can be equivalently defined as a pair
of sets $R,X$ with maps $s,t,e,c$ satisfying the various compatibilities.

If in addition every morphism in $C$ is an isomorphism, $C$ is
called a {\bf groupoid}. It is denoted by $R \rightrightarrows X$.
Then we have a map $i :R \rightarrow R$ sending a morphism to its
inverse.

Observe that if $C$ is a groupoid, then the image of
$j=(t, s) : R \rightarrow X \times X$ is an pre-equivalence relation
(where two objects are equivalent if they are isomorphic).

Let $s,t,e,c,R,X$ define a category, $C$.
If $A \rightarrow X$
is a map (of sets), then fiber product defines a category
$C|_A$, with objects $A$ and maps
$R|_A = R \times_{X \times X} A \times A$,
where a map between two objects $a,b$ in $A$ is just a map between
their images in $X$, with composition taking place
on $X$. If $C$ is a groupoid, so is $C|_A$. We refer to
$R|_A$ as the {\bf restriction} of $R$ to $A$.

Let $p: X \rightarrow Z$
be a  map (of sets) such that $ p \circ s = p \circ t$
and let $Z' \rightarrow Z$ be any map. Set $X' = X \times_Z Z'$
and $R'=R \times_Z Z'$. The two maps $s,t: R \rightarrow X$
induce $s',t':R' \rightarrow X'$, and it is immediate that
these define a groupoid. This is compatible with (1.7).
Note there is a natural injection
$R' \hookrightarrow R|_{X'}$, and that in this way
$(R',X')$ is a subcategory of $(R|_{X'},X')$. We refer to
$R'$ as the {\bf pullback}.

\definition{2.1 Definition} A {\bf groupoid space} is a quintuple of maps
$s,t, c, e, i$ (all of finite type) of algebraic spaces as above,
such that for all $T$, the quintuple
$s(T), t(T), c(T), e(T), i(T)$ define a groupoid
$R(T) \rightrightarrows X(T)$,
in a functorial way.

For a groupoid (space) $s,t : R \rightrightarrows X$, we use the following
convention throughout this paper: $j=(t,s), p_1=t, p_2=s.$
\enddefinition

 From now on we talk about algebraic spaces. For simplicity we
will abuse notation and call a groupoid space, a groupoid. When
we wish to return to the categorical setup, we will speak of a
groupoid of sets.

By a flat (resp \'etale, etc.)
groupoid, we mean one for which the maps $s,t$ are flat
(resp \'etale, etc.).

It is immediate from the categorical remarks
that a groupoid is an pre-equivalence relation, and
that the restriction or pullback (defined in the obvious way
generalizing the definition for sets) of
a groupoid
is again a groupoid. Note in particular that when $X$ is a geometric
point then $R$ is a group scheme, so in particular for any
geometric point $x \in X$, $R|_x$ is a group scheme, called
the {\bf stabilizer} at $x$. We call the group scheme
$S= j^{-1}(\Delta_X) \rightarrow X$ the stabilizer.

By a {\bf map } $g: (R',X') \rightarrow (R,X)$ {\bf between
groupoids} we mean a pair of maps $g:R' \rightarrow R$
and $g: X' \rightarrow X$ such that the induced map between
categories $(R'(T),X'(T)) \rightarrow (R(T),X(T))$ is a functor
for all $T$ (this of course can be equivalently expressed by
saying various diagrams commute).

We say a geometric point $x \in X$ is fixed by a pre-equivalence
relation, if $x$ is the only geometric point of its orbit
$t(s^{-1}(x))$.

It will be useful to know when a map of groupoids is a
pullback.

\definition{2.2  Definition} Let
$g:(R',X') \rightarrow (R,X)$ be a map of groupoids.

(1) We say that $g$ is a {\it square} when the commutative diagram
$$
\CD
R' @>>> X' \\
@V{g}VV        @V{g}VV \\
R  @>>>  X
\endCD
$$
is a fiber square, if we take for the horizontal maps either
both source maps, or both target maps.

(2) We say that $g$ is {\it  fixed point reflecting} if
for each geometric point $x' \in X'$ the induced map
of stabilizers $S_{x'} \rightarrow S_{g(x')}$ is a set
bijection.
\enddefinition

\remark{2.3 Remark} It is immediate that if $X \rightarrow Z$
is $R$-invariant, and $(R',X')$ is the pullback along any map
$Z' \rightarrow Z$, then $(R',X') \rightarrow (R,X)$ is square
and fixed point reflecting.
\endremark

Note that in (2.2.1), if the diagram is a fiber square for the source,
then the same follows for the target because $t=s \circ i$.

\definition{2.4} We will say that a \gc quotient $X \to Q$
for a groupoid $R \rightrightarrows X$ satisfies the {\bf descent
condition} if whenever $g:(R',X') \rightarrow (R,X)$ is an
\'etale square
fixed point reflecting map of groupoids such that
a \gc quotient $q' : X' \rightarrow Q'$ for
$X'/R'$ exists, then the induced map $Q' \rightarrow Q$ is \'etale
and  $X' \simeq X \times_Q Q'$.\enddefinition

\remark{2.5} One checks easily that a map of groupoids of
sets $g:(R',X') \rightarrow (R,X)$ is obtained by pullback from
the induced map $Q' \rightarrow Q$ iff $g$ is square and fixed
point reflecting. That is, the descent condition always holds for
groupoids of sets. \endremark

\remark{2.6 Remark} We note for convenience that for
any map $f: W \rightarrow X$, the restriction
$R|_W$ is described by the following diagram, in which every
square is a fiber square:
$$
\CD
R|_W @>>>  R \times_{(p_2,f)} W @>>> W \\
@VVV     @VVV   @V{f}VV \\
W\times_{(f,p_1)} R  @>>>   R @>{p_2}>> X \\
@VVV      @V{p_1}VV  @. @. \\
W   @>{f}>> X @. @.
\endCD
$$
\endremark

The following will alleviate separation anxiety:

\proclaim{2.7 Lemma-Definition}
Let $R \rightrightarrows X$ be a groupoid and let
$j : R \to X \times X$
be the induced map and $S$ the stabilizer group
scheme $j^{-1}(\Delta_X) \to X$ over X. Then $j$ is
separated iff $S \to X$ is separated.

In this case, we say that $R \rightrightarrows X$ has a
separated stabilizer.
\endproclaim

\demo{Proof} The only if part is obvious because $S \to X$
is the base change of $j$
by $\Delta_X \to X \times X$.

Assume that $S \to X$ is separated. Hence any section
in particular the identity section $e : X \to S$ is a
closed immersion (cf. \cite{Knudson71, II.3.11}).
Look at the commutative diagram:
$$
\CD
R \times_{X \times X} R @= R \times_{(s,s)} S \\
@A{\Delta}AA                  @AAA \\
R          @=                R \times_{(s,s)} e(X)
\endCD
$$
where the top identity is obtained by the isomorphism
$(r_1, r_2) \mapsto (r_1, i(r_1) \circ r_2)$.
Since $e(X) \to S$ is a closed embedding, so are the
vertical maps.
Thus $j$ is separated.
\qed
\enddemo

\proclaim{2.8  Corollary}
If $R \rightrightarrows X$ has a separated stabilizer,
so does $R|_W$ for any map $g : W \to X$.
\endproclaim

{\it We will assume throughout that every groupoid has
a separated stabilizer.}

\heading \S 3 Localizing the construction of quotients \endheading

\proclaim{3.1 Lemma} Let $R$ be an pre-equivalence relation.
Let $g:W \rightarrow X$ be a map, and $R|_W$
defined as in (2.6).

(1)
The canonical map $W/(R|_W) \rightarrow X/R$ of sheaves is
injective, and is an isomorphism if the composition
$$
p:W\times_{(g,p_1)} R  @>>>   R @>{p_2}>> X
$$
is surjective in the qff
topology.

(2) ({\bf Open}) Assume that $p_1, p_2 : R \to X$ are universally open,
that $p$ above is universally open, and as a map to its image,
is surjective in the qff topology. Suppose $X \rightarrow Z$
is a \gc quotient. Then the image, $V$, of the composition
$W \rightarrow X \rightarrow Z$ is open, and the induced map
$W \rightarrow V$ is a \gc quotient.

(3) ({\bf Surjection}) Assume that $p_1, p_2 : R \to X$
are universally open and
$p$ above is surjective in the qff topology.
Suppose $W/(R|_W) \rightarrow Z$ is a \gc quotient. By (1) there
is an induced map $X/R = W/(R|_W) \rightarrow Z$. It is a
\gc quotient.
\endproclaim

\demo{Proof}
Let $q_X$ and $q_W$ be the quotient maps (with domains $X$ and $W$).

If $Z' \rightarrow Z$ is a map,
then for $W',X',R'$ the fiber products of $W,X,R$ with $Z'$ over $Z$,
we have a fiber diagram
$$
\CD
W' \times_{g',p_1'} R' @>>> R' @>{p_2'}>> X' @>>> Z' \\
@VVV                        @VVV          @VVV    @VVV  \\
W\times_{g,p_1} R     @>>> R @>{p_2}>>   X   @>>> Z
\endCD
\tag{*}
$$
(here for a map $f$, $f'$ indicates the pullback).

(1) follows from the definitions of sheaves.

For (2): $V$ is open and $q_X$ is universally open, by (1.9).
Analogous remarks using $*$ show $q_X \circ g$ is universally open.

By definition the result holds if $W$ is an $R$-invariant open set,
thus we may replace $X$ by $U$. Then $W/(R|_W) = X/R$ so (1.8.G) and
(1.8.C) hold.

By $*$ the assumptions are preserved by any flat base extension
$Z' \rightarrow Z$. Thus (1.8.UC)
holds.

For (3):
By $*$ the assumptions are preserved by pullback along any flat
$Z' \rightarrow Z$, so it is enough to consider properties
(1.8.G), (1.8.C) and (1.8.US). (1.8.G) and (1.8.C) depend
only on the sheaf $W/(R|_W) = X/R$ and thus hold for $q_X$ iff they
hold for $q_W$. Since $q_W$ factors through $q_X$, if $q_W$
is a submersion, so is $q_X$. \qed \enddemo

\proclaim{3.2 Lemma} Assume that $p_1,p_2 : R \to X$ are
universally open. Suppose $\{U_i\}$ is a finite \'etale cover of
$X$ Suppose \gc quotients
$q_i: U_i/(R|_{U_i}) \rightarrow Q_i$ exist for all i. Then a
\gc quotient $q: X \rightarrow Q$ exists. \endproclaim

\demo{Proof} By (3.1.3), and induction, we may assume we have a
Zariski cover by two $R$-invariant open sets.
By (3.1.2), the \gc quotient of
$(U_1 \cap U_2)/(R|_{U_1 \cap U_2})$ exists as open set of the \gc
quotient of $U_i/(R|_{U_i})$ for $i=1,2$.
Since the categorical quotients are unique, these glue. The rest
is easy. \qed \enddemo

\proclaim{Corollary 3.2.1} Suppose $q: X/R \rightarrow Q$ is a map.
To check it is a \gc quotient we may work \'etale locally on $Q$.
\endproclaim
\demo{Proof} If $Q$ is \'etale locally the \gc quotient, then $q$
satisfies (1.8.C) by decent. Also $X$ has an \gc quotient by (3.2).
Thus since a categorical quotient is unique, $q$ is the \gc quotient.
\qed \enddemo

\proclaim{Lemma 3.3} Assume that $s, t : R \to X$ is flat and
$j$ is quasi-finite. Let $x \in X$ be a geometric closed point.
To prove that there is an \'etale neighborhood $U$ of $x \in X$
such that a \gc quotient exists for $U/(R|_U)$,
we may assume $R$ and $X$ are separated
schemes, $s, t : R \to X$ are quasi-finite and
flat and $x$ is fixed. To show the \gc quotient has any
additional property for which (3.1.2) and (3.1.3) hold,
we can make the same assumptions.
\endproclaim
\demo{Proof}
By (3.1.2) we can assume $X$ is a separated scheme. Since $R$ is
separated and quasi-finite over $X \times X$, $R$ is also a
separated scheme by \cite{Knudson71, II.6.16}.

\proclaim{Claim 1} We can assume $F = s^{-1}(x)$ is Cohen Macaulay
along $j^{-1}(x,x)$.
\endproclaim

\demo{Proof} Since $j$ is quasi-finite,
there is a geometric closed point $w \in F$ such that $F$
is Cohen Macaulay along $t^{-1}t(w)$.
(Keep cutting down $\overline{t(F)}$ with a general hypersurface
of high degree till one gets a 0-dimensional component $V$.
Then take $w \in F \cap t^{-1}(V)$.)
We have an isomorphism
$F @>{\circ i(w)}>> s^{-1}(t(w))$ sending $w$ to $e(t(w))$.
We can replace $x$ by $t(w)$, and thus the claim follows after
a Zariski shrinking around $x$.  \qed
\enddemo

\proclaim{Claim 2} We can assume $s,t$ are flat and quasi-finite.
\endproclaim
\demo{Proof} Let $x \in W \subset X$ be the closed subscheme
defined by $\xdim_x F$ elements of $m_x$, whose
intersection with $\ol{o_x}=\ol{t(s^{-1}(x))}$ is zero dimensional
at $x$
(lift any parameters of the maximal ideal of $\ol{o_x}$ at $x$).
Since $F$ is Cohen Macaulay along $j^{-1}(x,x)$, by \cite{Matsumura80}
$p: W \times_{(g,t)} R \overset {s} \to \rightarrow X$
is flat over  $x$, and  $s:R|_W \rightarrow W$
is flat and quasi-finite over $x$. After a Zariski shrinking, this
holds globally. The assertion on $t$ holds by $t=s \circ i$. By 3.1.3
we may restrict to $W$.
\qed \enddemo
Once $s,t$ are quasi-finite, we can assume $x$ is fixed after
a Zariski shrinking.
\qed \enddemo

\heading \S 4 Splitting \endheading

\definition{4.1 Definition} We say that a flat groupoid
$R \rightrightarrows X$ is {\bf split}
over a point
$x \in X$ if $R$ is a disjoint union of open and  closed subschemes
$R = P \du R_2$, with
$P$
a subgroupoid finite and flat, and
$j^{-1}(x,x) \subset P$.
\enddefinition

Our main goal in this section is to prove:

\proclaim{4.2 Proposition} Let $s,t:R \rightrightarrows X$ be a
quasi-finite flat groupoid of separated schemes.
Then every point $x \in X$ has an affine \'etale neighborhood
$(W,w)$ such that $R|_W$ is split over $w$.
\endproclaim

\proclaim{Lemma 4.3} Let $p:X \rightarrow Y$ be a map, and
$F,G \subset X$ two closed subschemes.
If $F$ is finite flat over $Y$ then there is unique closed
subscheme $i: I \rightarrow Y$ such that $f:T \rightarrow Y$
factors through $i$ iff $F_T$ is a subscheme of $G_T$.
\endproclaim
\demo{Proof} $I$ is obviously unique, so we can construct
it locally, and so can assume $p_*(\ring F.)$ is free. Then
there is a presentation
$$
K @>{h}>> p_*(\ring F.) \rightarrow p_*(\ring {F \cap G}.)
\rightarrow 0
$$
with $K$ free. $I$ is defined by the vanishing of $h$.
\qed \enddemo

Now let $s,t:R \rightrightarrows X$ be as in (4.2).
Note that $s : R \to X$ is quasi-affine
since it is separated and quasi-finite \cite{EGA, IV.18.12.12}.
Thus one can embed $s : R \to X$ into a projective scheme
of finite type over $X$. Thus the standard theory of Hilbert
schemes applies to our set up.

Let $g:H \rightarrow X$ be the
relative Hilbert scheme $\Hilb_{R/X}$ parametrizing
closed subschemes of $R$ which are proper flat over $X$
via $s$. Let $W \subset H$
be the closed subscheme of (4.3) of families containing
the identity section $e: X \rightarrow R$.
Let $s: P \rightarrow W$ be the universal family.

An $S$ point of $P$ consists of a pair $(F,a)$ of a map
$a:S \rightarrow R$, and a family
$F \subset R \times_{(s, s(a))} S$ flat over $S$, with $a$
factoring through $F$. Note that $g$ maps $(F,a)$ to
the $S$ point $[F]$ of $W$ representing
$F \subset R \times_{(s,s(a))} S$,
and $s$ sends $[F]$ to the $S$ point $s(a)$ of $X$.
The composition of component $R$ with $i(a)=a^{-1}$
gives an isomorphism $R\times_{(s,s(a))} S @>{\circ a^{-1}}>>
R \times_{(s,t(a))} S$. We write the image of $F$ under
this isomorphism as $F \circ a^{-1}$. Since $a$ factors through
$F$, the identity map $e$ factors through $F \circ a^{-1}$,
and $F \circ a^{-1} \subset R \times_{(s, t(a))} S$ defines
an $S$ point $[F \circ a^{-1}]$ of $W$. Thus we have a map
$t: P \rightarrow W$ given by
$t(F,a) = [F \circ a^{-1}]$, and $g([F \circ a^{-1}])=t(a)$.

\proclaim{4.4 Lemma} The pair $s,t: P \rightrightarrows W$ is
a finite flat subgroupoid of $R|_W$.
\endproclaim

\demo{Proof} In the big diagram (2.6) defining $R|_W$,
$P \subset R \times_{(s,g)} W$. Let $pr_R: P \rightarrow R$
be the induced map. Then
$g \circ t = t \circ pr_R$, and thus
$(F,a) \mapsto ([F \circ a^{-1}], a, [F])$
embeds $P$ in $R |_W$ (2.6).

Let $((F,a),(G,b))$ be an $S$ point of $P \times_{(s,t)} P$.
Then we have $F= G \circ b^{-1}$ and
$$
([F \circ a^{-1}], a, [F]) \circ ([G \circ b^{-1}], b, [G])
= ([G \circ (a \circ b)^{-1}], a \circ b, [G])
$$
by the composition law of $R |_W$.
Thus $P$ has the induced composition law
$$
(F,a) \circ (G,b) = (G, a \circ b)
 \quad \text{if} \quad F=G \circ b^{-1},
$$
the identity section is given by $e([F]) = (F,e)$,
and the inverse by
$i(F,a)=(F \circ a^{-1},a^{-1})$.
(Note $i(F,a) \in P$ by our definition of $W$,
this is the reason for restricting from $H$ to $W$.)
\qed \enddemo

Let $x \in X$ be an arbitrary point and let $P_x \subset R$ be
the open closed set of $R_x$ whose support is
$j^{-1}(x,x)$. Indeed if $M$ is the defining ideal of $x$
in $X$, then $P_x$ is defined in $R_x$ by
$(s^*M+t^*M)^n$ for $n \gg 0$.

Let $w=[P_x] \in W$, that is $P_w=P_x$.

\proclaim{4.5 Lemma}
$f : H \to X$ is \'etale at $w$, isomorphic on
residue fields $\xk(w) \simeq \xk(x)$ and
$H=W$ in a neighborhood of $w$.
\endproclaim
\demo{Proof}
The construction of $P_w \subset R \times_X \xk(x)$
shows $\xk(w) \simeq \xk(x)$.
It remains to see
$$
Hom_{(X,x)}((S,y),(W,w)) = Hom_{(X,x)}((S,y),(H,w))
=(\text{one point set)}
$$
for any artin schemes $(S,y)$ over $(X,x)$.
Since $R_S$ contains an open closed subscheme $P'$ such that
$P' \times_S y=P_x \times_x y$, we see $P' \supset e(S)$ and
these are obvious.
\qed
\enddemo

\proclaim{4.6 Corollary} $R|_W = P \cup R_2$, for $R_2$ a
closed subscheme
with $R_2 \cap P \cap s^{-1}(w) = \emptyset$.
\endproclaim
\demo{Proof} Let $R' = R|_W$, and let $y$ be an arbitrary point of
the finite set $P_w$. Then
$P_w = R'_w$ as schemes at $y$.
Let $I \subset \ring R'.$ be the defining ideal of $P$ in $R'$.
Then $P_w = R'_w$ at $y$ implies $I \otimes_{\ring W.} \xk(w)=0$
at $y$ by the flatness of $P$ over $W$. Thus
$I \otimes_{\ring R'.} \xk(y)=0$. Hence in a neighborhood of
$y$, we have $I=0$ (that is $P=R'$) by Nakayama's Lemma.
The result follows. \qed \enddemo

Now to prove (4.2), we will shrink $W$ so that $R_2$ and
$P$ become disjoint.  In order to preserve
the finiteness of $P$, we want to shrink using $P$ invariant open
sets. For this a general construction will be useful:

\subhead{ A geometric construction of equivalence classes} \endsubhead

\proclaim{4.7 Lemma} Let $A \rightarrow B$ be a local map
of local Noetherian rings, with
$B/(m_A \cdot B)$ of finite dimension $k$ over $A /m_A$.
Then $m_B^k \subset m_A \cdot B$.
\endproclaim

\demo{Proof} We can replace $B$ by $B/(m_A \cdot B)$, and so
assume $A$ is a field. Now the result follows by Nakayama's
Lemma. \qed \enddemo

Let $P \rightrightarrows X$ be a finite flat groupoid, $x \in X$
a point fixed by $P$, and
$k$ the degree of $s$ over $x$.

We consider the functor on $X$-schemes, whose $(\alpha : T \to X)$-point
consists of $k$ sections $g_i:P \times_{(s,\alpha)} T \to T$ of
$P$ whose {\it scheme theoretic union}
in $T \times R$ (the closed subscheme defined by the product of
the defining ideals of the $k$ sections $g_i$)
contains the pullback $\Gamma_T$, where $\Gamma \subset X \times P$
is the graph of $s$.
By 4.3, this functor is represented by a closed subscheme $I$
of the $k$-fold fiber product of $s : P \to X$.

A geometric point of $I$ consists of a $k$-tuple of points in
a fiber of $s$, whose scheme theoretic union contains the
fiber. Let $p:I \rightarrow X$ be the map induced by $s$.
By 4.7, $p$ is surjective on a neighborhood of $x$.
Let $\pi_i$ be the map
$$
\pi_i : I \subset P^{\times_X k}
( = \oversetbrace \text{$k$ times} \to {P \times_X \cdots \times_X P})
@>{pr_i}>> P @>{t}>> X.
$$
For a point $\alpha \in I$, note that the set
$\bigcup_{i}^k \pi_i(\alpha) \subset X$ is the full $P$-equivalence
class $[p(\alpha)]$.

\proclaim{4.8 Lemma} With the above notation and assumptions,
let $\xi$ be a geometric point at $x$ and  $U$ an
open set $U \supset [\xi]$ such that $s$ is of constant
degree $k$ over $U$. Let $V = s(t^{-1}(U^c))^c$. Then
$V$ is an $P$-invariant open neighborhood of $x$,
$V \subset U$, and $V = \{\eta \in X | [\eta] \subset U\}$.
If $U$ is affine then $V$
is affine. In particular, $x$ has
a base of $P$-invariant open affine neighborhoods.
\endproclaim
\demo{Proof} Since the projections are closed, everything before
the final assertion is clear.

Now assume $U$ is affine. Set
$$
J = \bigcap_{i=1}^{k} \pi_i^{-1}(U).
$$

Note (as a set) $J =\{\alpha | [s(\alpha)] \subset U \}$.
Since $p$ is surjective, $J = p^{-1}(V)$, and hence
$p:J \rightarrow V$ is finite. Since $\pi_i$ is finite, and $I$ is
separated, $J$ is affine. Thus so is $V$.
\qed \enddemo

\demo{Proof of (4.2)}
After a Zariski shrinking, we can assume $x$ is fixed
by $R$. We have by (4.6),
$R|_W = P \cup R_2$, with $w$ outside the
image (under either projection) of $R_2$.
This is preserved by restriction to $P$-invariant  open sets.
Since $R_2 \cap P \subset P$ is closed, and does
not meet the fiber, we can assume
by (4.8) and (4.5), that $R_2$ and $P$ are disjoint, and
$W \rightarrow X$ is \'etale. \qed
\enddemo

\heading \S 5 quotients for finite flat groupoids with affine base
\endheading

\proclaim{5.1 Proposition (Finite-Over-Affine case)}
Let $A$ be a Noetherian ring and $B$ an $A$-algebra
of finite type. Let $R$ be an affine groupoid finite and free over
$\Spec B$. Then $B^R$
is an $A$-algebra of finite type and
$$
q: X=\Spec B \rightarrow Q=\Spec(B^R)
$$
is finite and the \gc quotient of $X/R$.
\endproclaim

For some large number $n$, let $x_1=1, x_2, \cdots, x_n \in B$
be a set of generators as an $A$-algebra.
Let $\xi_1, \cdots, \xi_n$ be indeterminates over $B$
and by the flat base change $A \to
A[\xi_1, \cdots, \xi_n]$ we may treat $\xi$'s
as $R$-invariant indeterminates. Let $\phi=\sum_i x_i \xi_i$.

\proclaim{5.2 Lemma} We have $Nm_{t} (s^*\phi) \in B^R[\xi]$,
where $Nm_t$ is the norm for the finite free morphism $t : R \to X$.
\endproclaim
\demo{Proof of (5.2)}
It is enough to prove $t^*Nm_{t} (s^*\phi)=s^*Nm_{t} (s^*\phi)$.
By the definition of norm, we have
$$
s^*Nm_{t} (s^*\phi)=Nm_{pr_1}({pr_2}^*s^*\phi),
\tag{5.2.1}
$$
where $pr_1, pr_2$ are defined in the following diagram.
$$
\CD
R \times_{(s,t)} R @>{pr_2}>> R @>s>> X \\
@V{pr_1}VV @V{t}VV @. \\
R @>s>> X @. \\
\endCD
$$
Let
$$
\tau : R \times_{(s,t)} R \simeq R \times_{(t,t)} R
$$
be a morphism defined by
$(r_1,r_2) \mapsto (r_1,r_1 \circ r_2)$, which is an
$R$-isomorphism via the first projection.
Then we have a commutative
diagram.
$$
\CD
R @<{pr_1}<< R \times_{(s,t)} R @>{pr_2}>> R @>s>> X \\
@|          @V{\tau}VV                      @.    @| \\
R @<{pr'_1}<< R \times_{(t,t)} R @>{pr'_2}>> R @>s>> X \\
\endCD
$$
Similarly to (5.2.1), we see
$t^*Nm_{t} (s^*\phi)=Nm_{pr'_1}({pr'_2}^*s^*\phi)$. Then
starting with $\phi$ on $X$ on the right, we can send it along the
diagram either by pull back or by norm. On the left, we get the
equality
$$
s^*Nm_{t} (s^*\phi)=Nm_{pr_1}({pr_2}^*s^*\phi)
=Nm_{pr'_1}({pr'_2}^*s^*\phi)=t^*Nm_{t} (s^*\phi),
$$
which proves the $R$-invariance as required. \qed
\enddemo

\proclaim{5.3 Lemma} $B^R$ is an $A$-algebra of finite type and
$B$ is a finite $B^R$-module. In particular $q$ is finite and surjective.
The formation of $Q$ commutes with any flat base extension
$Q' \rightarrow Q$ with $Q'$ affine.
\endproclaim

\demo{Proof of (5.3)}
Let $F(\xi)=Nm_t(s^*\phi) \in B^R[\xi]$ by (5.2).
Let $C$ be the $A$-subalgebra of $B^R$ generated by all the
coefficients of $F$. We note that
$$
F(x_i,0,\cdots,0,{\overset{i\text{-th}} \to -1},0,\cdots,0)=0
\qquad \forall i \ge 2
$$
because $t : R \to X$ has a section $e$ such that
$s \circ e= t \circ e = id$. Since $F(\xi)$ is monic in $\xi_1$,
each $x_i$ is integral over $C$. Since $x_i$ generate $B$ as an
$A$-algebra, $B$ is finite over $C$.
Thus $B^R$ is also a finite $C$-module and hence the
lemma is proved. \qed
\enddemo

\proclaim{5.4 Lemma}
$q$ satisfies the condition (1.8.G).
\endproclaim

\demo{Proof of (5.4)}
Since $q$ is finite dominating, the map in (1.8.G) is surjective.
It is enough to prove the injectivity.
Let us look at $F(\xi)$ in the proof of (5.3)
as $F_{t}(\xi)$ parametered by $t \in \Spec B^R$.
Let $a \in X$ be a geometric point. Then by the definition of
norm, we have a polynomial identity
$$
F_{q(a)}(\xi) = \prod_{b \in t^{-1}(a)}
(\xi_1+\sum_{i=2}^n x_i(s(b)) \xi_i),
$$
where $b$ is chosen with multiplicity. So if $q(a)=q(a')$,
then there exists
$b \in t^{-1}(a)$ such that $x_i(s(b))=x_i(s(a'))$ for all $i$.
Since $x_i$ generate the $A$-algebra $B$, we have $a'=b$.
\qed \enddemo

\demo{Proof of (5.1)}
By 3.2.1 and 5.3 we can work locally in the \'etale topology on $Q$.

For the condition (1.8.UC), let $g : X \rightarrow Z$
be an $R$-invariant map. It is enough to
prove that for any geometric point $q_0 \in  Q$ there is an \'etale
neighborhood $Q'$ of $q_0$ such that $X \times_Q Q' \to Z$ factors
through $Q'$. Passing
to the strict henselization, we may assume $(Q,q_0)$ is strictly
henselian. Then
$X$, since it is finite over $Q$, is a disjoint union of strictly
henselian local schemes, and $R$ acts transitively on the components
by (5.4).
By (3.1.3), we may drop all but one component without affecting $Q$,
or the
existence of a factorization, and so assume $X$ is strictly
henselian $(X,x)$. Then $X \to Z$ factors through a strict
henselization of
$(Z,g(x))$ which is $R$-invariant since $R$ acts trivially on the
residue field of $X$. Thus $g$ factors through $Q$
by construction of $Q$.

$q$ is a universal submersion since it is a finite surjection.
\qed \enddemo

\heading \S 6 Auxiliary Results \endheading

Let $g:(R',X') \rightarrow (R,X)$ be a map of groupoids. Let
$T \rightarrow X$ be a map and $T' = T \times_X X'$. There is
then an induced map of groupoids
$g_T:(R'|_{T'},T') \rightarrow (R|_T,T)$.

\proclaim{6.1 Lemma} If $g$ is square, so is $g_T$. \endproclaim
\demo{Proof} Passing to $\Hom$ we can assume we are working with
groupoids of sets. One simply considers the natural map
$R'|_{T'} @>{g_T \times s}>> R|_T \times_T T'$ and checks
it is a bijection.
\qed \enddemo

\proclaim{6.2 Remark} Let $f : X \to Y$ be a map of
algebraic spaces.
Then

(1) $g : X \times_Y X \to X \times X$ is separated.

(2) If $X$ is separated, then $f$ is separated.
\endproclaim
\demo{Proof} The map $g$ is a monomorphism
(of sheaves) by the
definition of $X \times_Y X$. Hence $g$ is separated.
Assume that $X$ is separated.
Since the composition $X \to X \times_Y X \to X \times X$ is
a closed immersion, so is $X \to X \times_Y X$
\cite{Knudson71, I.1.21}.
\qed
\enddemo

\proclaim{6.3 Lemma} Let $R \rightrightarrows X$ be a
finite flat groupoid with separated $X$.
Let $q:X \rightarrow Q$ be a \gc quotient. Then $q$
is finite, it is the \gc quotient, and
also satisfies the descent condition (2.4).  \endproclaim
\demo{Proof} Since $X$ is separated, so is $q$ (6.2).
By (3.2.1) we can assume $Q$ is an affine scheme by a base change.
Since $q$ is separated and quasi-finite by (1.8.G), $X$ and
$R$ are separated schemes \cite{Knudson71, I.I6.16}.

After shrinking $Q$ we can assume
$q$ is quasi-projective \cite{EGA, IV.18.12.12}:
We can complete $q$ to $X \subset X' @>{q'}>> Q$
with $X \subset X'$ open and $X' @>{q'}>> Q$ projective.

Let $F \subset X$ be the fiber over a geometric point $p \in Q$.
Then $F \subset X'$ is closed. Let $Z = X' \setminus X$, and
let $H \subset X'$ be a general hypersurface such that $Z \subset H$.
Then $H \cap F = \emptyset$. Thus $F\subset H^c=U \subset X$ is an
open affine set. By 4.8 we can assume $U$ is $R$-invariant, and
so can assume $X$ is affine. Thus $q$ is finite, and so it
is a universal submersion. The additional properties follow from
5.1.

Now
let $g:(R',X') \rightarrow (R,X)$ be as in the definition 2.4.
We can check the descent condition locally
in the \'etale topology on $Q$. Thus we may assume $Q=(Q,z)$
is strictly henselian and $x' \in X'$ any point lying over $z$.
Then $X$ is a disjoint
union of strictly henselian schemes.  We can obviously prove
$X ' = X \times_Q Q'$ one connected component of $X$ at a time.
Note that if $W$ is
a connected component of $X$, then $R|_W \rightrightarrows W$
is again finite. By (3.2) replacing $X$ by a connected component
containing $g(x')$ does not affect $Q' \rightarrow Q$, and by
(6.1) the assumptions are preserved.
So we may
assume $X$ is strictly henselian. Then $X'$ is a disjoint union
$\du X_i \du Y_1$ with each $X_i$ isomorphic under $g$ to $X$,
and $z \not \in g(Y_1)$. Let $i$ be such that
$\{x'\}=g^{-1}(x) \cap X_i$.
Since $R' \rightrightarrows X'$ is finite flat,
$t'({s'}^{-1} X_i)$ is an open closed
set of $X'$ and finite over $X$. Thus $t'({s'}^{-1} X_i)$
is a union of some $X_j$'s.
The fixed point reflecting
condition $S_x = S_{x'}$ means that
$g^{-1}(x) \cap t'({s'}^{-1} X_i)=\{x'\}$.
Thus $t'({s'}^{-1} X_i)=X_i$, and $X_i$ is $R'$-invariant. Thus
$(R'|_{X_i},X_i) \rightarrow (R,X)$ is square,
and we can assume $X'=X_i$.
Then we have $(R',X') \simeq (R,X)$.
 \qed
\enddemo

\proclaim{6.4 Lemma} Let $s,t : R \to X$ be a
groupoid with $j$ proper. Then
if $q : X \to Q$ is a \gc quotient for $X/R$, then $Q$ is separated.
\endproclaim
\demo{Proof} We can replace $X$ with its \'etale affine cover
(3.1.2). The proof of (2.9) of \cite{Koll\'ar95} extends
to groupoids without change.
\qed \enddemo

\proclaim{6.5 Fine Quotients Lemma}
Let $R \rightrightarrows X$ be a finite flat groupoid
of separated schemes such that $j : R \to X \times X$ is a closed
embedding. Then $X/R$ is represented by a finite flat
morphism
$q : X \to Q$ for the qff topology.
The map $q$ is flat, the construction
commutes with pullback along any $Q' \rightarrow Q$, and
$q$ is a \gc quotient.
\endproclaim
\demo{Proof}
Since $R \subset X \times X$ is finite and flat over $X$,
$t:R \rightarrow X$
defines an $R$-invariant map $q: X \rightarrow \Hilb_X$.
Let $i:Z \subset \Hilb_X \times X$ be the universal family. Then
$R = X \times_{\Hilb} Z$. By the universal property of $Z$,
the section $e$ of $t$ induces a map
$\gamma: X \rightarrow Z$ such that $i \circ \gamma = (q,id)$.
thus $\gamma$ is a closed embedding. In particular, since
$\pi:Z \rightarrow \Hilb_X$ is finite, $q = \pi \circ \gamma$
is finite. Let $Q \subset \Hilb_X$ be the scheme
theoretic image of $q$.
We abuse
notation, and write for $Z$ its restriction to $X \times Q$,
whence $Z \times_Q X = R \subset X \times X$.

Let $q:R \rightarrow Z$ be the projection.

We now show $\gamma: X \hookrightarrow Z$ is an isomorphism.
Since $\gamma$ is a closed embedding,
$X \times_Q X \subscheme X \times X$
factors as
$$
X \times_Q X = X \times_Q Z \times_Z \gamma(X) =
R \times_Z \gamma(X) \subscheme R \subscheme X \times X.
$$
Since $q$ is $R$-invariant, $R \subset X \times_Q X$. Thus
$R = R \times_Z \gamma(X)$.

Since $q$ is affine, by
the definition of the scheme theoretic image,
$\ring Q. \rightarrow q_*(\ring X.)$
is injective. Since $\pi$ is flat,
$\ring Z. \rightarrow q_*(\ring R.)$
is injective. Thus $R = R \times_Z \gamma(X)$ implies $\gamma$
is an isomorphism.

Thus we have that $q$ is flat and $R = X \times_Q X$. Note these
two properties are preserved by pullback along any
$Q' \rightarrow Q$.
They also imply $Q = X/R$ as sheaves in the qff topology.
Now the rest follows easily. \qed
\enddemo

\heading \S 7 The Boot Strap Theorem \endheading
We begin with a number of categorical remarks, which
will guide our construction when we move to algebraic
spaces. In addition, for \gc quotients, these categorical
remarks will describe what happens on $T$-valued points
for arbitrary schemes $T$.

Suppose
$s,t:R \rightarrow X$ defines a groupoid of sets.

Let $S \subset R$ be the subgroupoid of automorphism.
Let $P \subset R$ be a subgroupoid, and $S_P = S \cap P$.

\definition{7.1 Definition} $P$ is {\bf normal} if the
map
$$
S_P \times_{(s,t)} R \rightarrow S
$$
given by $(s,r) \rightarrow i(r) \circ s \circ r$ factors
through $S_P$.
\enddefinition

(7.2) $P$ acts on $R$ by either pre-composition, or
post-composition, or both. Thus we have groupoids (of sets):
\roster
\item $P \rightrightarrows X$
\item $P \times_{(s,t)} R \rightrightarrows R $
\item $R \times_{(s,t)} P \rightrightarrows R $
\item $P \times_{(s,t)} R \times_{(s,t)} P \rightrightarrows R$
\endroster

(7.3) We indicate their quotients respectively by $X',R^t,R^s,R''$.
Note there are natural maps $s: R^t \rightarrow X$,
$t: R^s \rightarrow X$, and $s'',t'': R'' \rightarrow X'$.
and groupoids
\roster
\item $P \times_{(s,t)} R^s \rightrightarrows R^s $
\item $R^t \times_{(s,t)} P  \rightrightarrows R^t .$
\endroster

It is elementary to check that $R''$ is the quotient of either
7.3.1 or 7.3.2.

The following is easy (if a bit tedious) to check:

\proclaim{7.4 Lemma} $s'',t''$ define a groupoid such that
the map $(R,X) \rightarrow (R'',X')$ is a map of groupoids
iff $P$ is normal.
\endproclaim

Pre-composition on the first factor, and post composition
composed with the inverse on the second factor give a groupoid
$$
M=P \times_{(t,s \circ pr_1)} (R^t \times_{(s,t)} R^s)
\rightrightarrows R^t \times_{(s,t)} R^s \tag{7.5}
$$
(informally, the action is
$(p,a,b) \rightarrow (a \circ p, i(p) \circ b)$).
There are maps of groupoids
$$
\align
g:(R^t \times_{(s,t)} P, R^t) &\rightarrow (P,X) \\
h:(M, R^t \times_{(s,t)} R^s) &\rightarrow (R^t \times_{(s,t)} P, R^t)
\endalign
$$
which are  square by construction.

\proclaim{7.6 Lemma} If $P$ is normal then $g$ and $h$  are fixed point
reflecting.
\endproclaim
\demo{Proof}
First for any $P$, note that for $\beta \in R$, and
$p \in P$, $[\beta \circ p] = [\beta]$ in $R^t$
iff $\beta \circ p \circ i(\beta) \in P$.
Thus if $P$ is normal, the stabilizer of $[\beta]$ is the fiber
of $S_P$ over $s(\beta)$. Of course similar result apply to
$R^s$. Now the lemma follows easily.
\qed \enddemo

Let $Z = R^t \times_{(s,t)} R^s / M$.
There is a natural $M$ invariant surjection
$$
R^t \times_{(s,t)} R^s \rightarrow R'' \times_{(s'',t'')} R''
$$
and thus a surjection $p: Z \rightarrow R'' \times_{(s'',t'')} R''$.
\proclaim{7.7 Lemma} If $P$ is normal, $p$ is an isomorphism. If
$S \subset P$ then $R'' \rightarrow X' \times X'$ is injective.
\endproclaim
\demo{Proof}
This follows easily from the stabilizer remarks in the proof of
7.6. \qed \enddemo

{\bf Now we work with separated schemes:}

\proclaim{7.8 Boot Strap Theorem} We follow the above notation.
Assume $R \rightrightarrows X$, is a quasi-finite flat groupoid
of separated schemes, and $P \subset R$
a closed and open subgroupoid which is finite (and necessarily flat)
over $X$. Assume \gc quotients exist in (7.2.1), (7.2.4), (7.3.1)
and (7.3.2) (Note quotients exist for (7.2.2) and (7.2.3) by (6.5)).
\roster
\item If $S \subset P$, then $R'' \rightrightarrows X'$ is
an \'etale equivalence relation and the algebraic space
$X'/R''$ is the \gc quotient for $X/R$.
\item If $P$ is normal, $R'' \rightrightarrows X'$ is an
an \'etale groupoid such that \gc quotient for $X/R$ exists iff
one for $X'/R''$ does, and if it does then they are isomorphic.
I.e. $X/R$ and $X'/R''$ define the same \gc quotient problem
in a neighborhood of $x$.
\endroster
\endproclaim
\demo{Proof}
Note the groupoids (7.2.2) and (7.2.3)  are free, and thus fine
quotients exist by (6.5). The map
$$
R \times_{(s,t)} R \rightarrow R^t \times_{(s,t)} R^s
$$
is also a fine quotient, as it is obtained from the fine
quotient $R \times R \rightarrow R^t \times R^s$
by pullback. In particular these quotients are flat and
commute with any base extension.

We show next $s: R^t \rightarrow X$ is \'etale. It is flat since
$R \rightarrow R^t$ is flat.
$e$ gives a section $\ol{e}$ of $s$, and $\ol{e}(X)$ is a
connected component
since $P$ is a connected component of $R$. Thus $s$ is unramified along
$\ol{e}$. Translating by the inverse gives an automorphism of the
fiber carrying a prescribed point to $\ol{e}$, thus $s$ is unramified,
and hence \'etale.

Now $g$ and $h$ are \'etale, square, and fixed point reflecting by (7.6)
(the only thing left to check is fixed point reflecting, which
is a set theoretic question, and thus follows from the categorical
case since the quotients are geometric).

Now by (6.3)
$s''$, $t''$, $\ol{h}: Z \rightarrow R''$ are
\'etale.  The universal property of $Z$ induces a commutative diagram
$$
\CD
Z @>{p}>>  R'' \times_{(s'',t'')} R'' \\
@V{\ol{h}}VV  @V{pr_1}VV \\
R'' @= R'' .
\endCD
$$
$p$ is a bijection on geometric points by (7.7).
Since both sides are  \'etale over $R''$, the map is an isomorphism.
The universal properties of both maps
$$
R\times_{(s,t)} R \rightarrow R^t \times_{(s,t)} R^s
\rightarrow Z
$$
together with the composition of $R$ induce a map
$c: R'' \times_{(s'',t'')} R'' \rightarrow R''$. Such that
the diagram
$$
\CD
R \times_{(s,t)} R @>{c}>> R \\
@VVV             @VVV         \\
R'' \times_{(s'',t'')} R'' @>{c}>> R ''
\endCD
$$
commutes. Similarly we have induced commutative diagrams
$$
\CD
R @>{i}>> R \\
@VVV      @VVV \\
R'' @>{i}>> R''
\endCD
$$
and
$$
\CD
X @>{e}>> R \\
@VVV      @VVV \\
X' @>{e}>> R''.
\endCD
$$
To check these define a groupoid, we need to check various
diagrams commute. In most cases this follows from the corresponding
diagram of $(R,X)$ by the universal property. For example to
see
$$
\CD
R'' \times_{(s'',t'')} R'' @>{c}>> R '' \\
@V{pr_1}VV                 @V{t''}VV \\
R''  @>{t''}>>             X'
\endCD
$$
commutes, one compares this with the commutative diagram
$$
\CD
R \times_{(s,t)} R @>{c}>> R  \\
@V{pr_1}VV                 @V{t}VV \\
R  @>{t}>>             X
\endCD
$$
and uses that
$$
\align
R \times_{(s,t)} R & \rightarrow R'' \times_{(s'',t'')} R'' \\
R & \rightarrow R'' \\
X & \rightarrow X \\
\endalign
$$
are all categorical surjections, by the universal properties.
Note in particular that since $t''$ and $s''$ are \'etale, this
gives that $c$ is \'etale.  The only diagram which must be
treated differently is for the associativity of composition. Since
$c$ is \'etale, this
amounts to checking equality between \'etale maps. This is a
set theoretic question, and so follows from (7.4). One could also express
the three fold fiber products of $R''$ as geometric quotients and then
apply the universal property.

Now suppose $S \subset P$. $R'' \rightarrow X' \times X'$ is unramified,
since the projections are \'etale. It is injective on geometric points
by (7.4). Thus it is a monomorphism. Now let $Q = X'/R''$. Since
$q:X \rightarrow Q$ is the composition of two \gc quotient
maps, the result follows easily.
\qed \enddemo

\heading \S 8 Proof of Theorem 1.1 \endheading

\proclaim{8.1 Lemma} Let $R \rightrightarrows X$ be a flat quasi-finite
groupoid with finite $j$. Then a \gc quotient exists.
\endproclaim
\demo{Proof} Note the assumptions are preserved by any flat
pullback. By (3.2) its enough to construct the quotient locally
in a neighborhood of $x \in X$.
 Thus we can assume $x$ is fixed
by $R$. By (4.1) we can assume $R$ is split over $x$, and that
$X$, and hence $R$ (since $j$ is finite)
are affine schemes. Then each of the necessary \gc
quotients in (7.8) exists by (5.1) \qed \enddemo

\demo{Proof of Theorem 1} Since the stabilizer (2.7) is finite,
$j$ is quasi-finite. We can assume we are in the situation
of (3.3). We reduce to $R$ split as in the proof of (8.1). By
(8.1) each of the necessary quotients in (7.8) exists. \qed \enddemo

\heading \S 9 Relations to GIT \endheading
Throughout this section we assume $L$ is the spectrum of a field.
We follow the notation of \cite{MumfordFogarty82}.

Note first that for any map $p:S\rightarrow X$ there is a
unique
maximal open subset of $X$ over which $p$ is proper. When $p$
is the stabilizer of a groupoid, we call this open set
$X(PropStab)$.

\proclaim{9.1 Proposition} Let $G$ be a reductive group acting on
a scheme $X$. Then $U= X^s_0(Pre) \subset X(PropStab)$, and the quotients
of $U$ given by \cite{MumfordFogarty82,1.9} is isomorphic to the \gc
quotient of $U/G$.
\endproclaim
\demo{Proof} Let $q: U \rightarrow Y$ be the GIT quotient.
To check finiteness of the stabilizer, we can replace
$Y$ by any etale cover, and so can assume $Y$
and thus $U$ are affine (note $q$ is affine), and in particular
separated. Then by \cite{MumfordFogarty82,0.8} the
action of $G$ is proper.

Now let $g: U \rightarrow Q$ be the GC quotient of (1.1). $Q$
is universal among algebraic spaces, while $Y$ is universal
among schemes. We have an induced map $h: Q \rightarrow Y$.
We can check $h$ is an isomorphism locally on $Y$, and
so can assume as above that $G$ acts properly.
Then $Q$ is separated, and $h$ is set theoretically one to
one. Thus by \cite{Knudson71, I.I6.16} $Q$ is a scheme, and
so $h$ is an isomorphism by the universality of $Y$. \qed \enddemo

\remark{9.2 Remark} A reductive group acting with quasi-finite
stabilizer does not in general have finite stabilizer. For example if
$G = Aut(P^1)$ acting on the $n$-th symmetric product $X$ of $P^1$
for $n \ge 4$ and $U$ is the open set of $X$ where the stabilizers
are quasi-finite, then $U \not \subset X(PropStab)$. \endremark

\remark{9.3 Remark} In (1.1), assume that
$j$ is finite, $s,t$ are affine, and
$L$ is the spectrum of a field.
Then the argument of
\cite{MumfordFogarty82,0.7}, which extends without change
to groupoids, shows that the quotient
map is affine. Mumford also remarks on pg. 16 that
(0.7) holds for arbitrary base. \endremark

\bigpagebreak
Mathematics Department,
University of Texas at Austin
Austin, Texas 78712 USA

email: keel\@mail.ma.utexas.edu

\medpagebreak

Research Institute for Mathematical Sciences,
Kyoto University,
Kyoto, 606 Japan

email : mori\@kurims.kyoto-u.ac.jp
\end{document}